\documentclass[12pt,a4paper,oneside,notitlepage]{article}
\usepackage[utf8]{inputenc}            
\usepackage[T1]{fontenc}

\usepackage{graphicx} 
\usepackage{geometry}
\usepackage{amsmath}             
\usepackage{amssymb}             
\usepackage{amsfonts}
                                 
\usepackage{latexsym}            
\usepackage{amsthm}             
\usepackage{cases}                       
\usepackage{mathrsfs}
\usepackage{xcolor}
\usepackage{lineno}
\usepackage{tikz}
\usepackage{hyperref}
                           
\newcommand{\field}[1]{\mathbb{#1}}        
\newcommand{\R}{\field{R}}

\newcommand{\Z}{\field{Z}}

\newcommand{\rmi}{\mathrm{i}}
\newcommand{\e}{\mathrm{e}}

\renewcommand{\d}{\mathrm{d}}

\newcommand{\pt}{p_t^\delta}
\newcommand{\xad}{x_\alpha^\delta}
\newcommand{\xd}{x^\dagger}
\newcommand{\pd}{p^\dagger}
\newcommand{\D}{\mathbf{D}}

\theoremstyle{plain}
\newtheorem{theorem}{Theorem}[section]

\newtheorem{corollary}[theorem]{Corollary}

\newtheorem{proposition}[theorem]{Proposition}
\theoremstyle{definition}
\newtheorem{definition}[theorem]{Definition}

\theoremstyle{remark}
\newtheorem{remark}[theorem]{Remark}

\title{Fast reconstruction approaches for photoacoustic tomography with smoothing Sobolev/Mat{\'e}rn priors}

\author{Jaakko Kultima\footnote{Johann Radon Institute for Computational and Applied Mathematics (RICAM), Linz, Altenbergerstraße 69, A-4040 Linz, Austria,
(jaakko.kultima@ricam.oeaw.ac.at)},
Ronny Ramlau\footnote{Johann Radon Institute for Computational and Applied Mathematics (RICAM), and Industrial Mathematics Institute, Johannes Kepler University, Linz, Altenbergerstraße 69, A-4040 Linz, Austria,
(ronny.ramlau@ricam.oeaw.ac.at)},
Teemu Sahlström\footnote{Department of Technical Physics, University of Eastern Finland, P.O. Box 1627, 70211 Kuopio, Finland (teemu.sahlstrom@uef.fi)},
Tanja Tarvainen\footnote{Department of Technical Physics, University of Eastern Finland, P.O. Box 1627, 70211 Kuopio, Finland, (tanja.tarvainen@uef.fi)}}

\date{\today}

\begin{document}
\maketitle

\begin{abstract}
    In photoacoustic tomography (PAT), the computation of the initial pressure distribution within an object from its time-dependent boundary measurements over time is considered. This problem can be approached from two well-established points of view: deterministically using regularisation methods, or stochastically using the Bayesian framework. Both approaches frequently require the solution of a variational problem. In the paper we elaborate the connection between these approaches by establishing the equivalence between a smoothing Mat{\'e}rn class of covariance operators and Sobolev embedding operator $E_s: H^s \hookrightarrow L^2$. We further discuss the use of a Wavelet-based implementation of the adjoint operator $E_s^*$ which also allows for efficient evaluations for certain Mat{\'e}rn covariance operators, leading to efficient implementations both in terms of computational effort as well as memory requirements. The proposed methods are validated with reconstructions for the photoacoustic problem.
\end{abstract}

\section{Introduction}
Photoacoustic tomography (PAT) is a hybrid imaging modality exploiting the photoacoustic effect. A pressure distribution is generated by a short light pulse, which then relaxes as an acoustic pressure wave. The pressure distribution is caused by the molecular distributions of the target with varying light absorption and thermal expansion properties \cite{li2009,beard2011,cox2012a,tarvainen2024quantitative}. PAT is able to provide 3D images of biological tissues with high-resolution which makes it very appealing for various applications \cite{li2009,beard2011}.

For the inverse problem of PAT, the generation of the pressure is assumed to be  instantaneous, and the initial pressure is estimated from the pressure wave that is measured on the boundary of the target using ultrasound sensors. 
The inverse problem is well-posed if the measurement domain surrounds the whole imaged target and the time-series is measured long enough, i.e., data exists for all times where the pressure wave is non-zero \cite{xu2004b}. 
Depending on the application, the qualitative information provided by the reconstruction of the pressure distribution, such as shapes, sizes and locations of inclusions within the target, may be sufficient. The problem of recovering the precise distributions of optical parameters in the tissue is referred to as quantitative PAT \cite{tarvainen2024quantitative,suhonen2023}. While beyond the scope of the present paper, we note that PAT can be viewed as the first of two steps when considering the quantitative PAT problems, and can thus benefit from efficient algorithms for the reconstruction of the pressure distribution. 

If a sufficiently large part of the boundary is covered by ultrasound sensors, reconstructions with sufficient accuracy can be obtained using backprojection-type algorithms. 
In addition to the 'universal backprojection' \cite{xu2005}, solutions to special geometries  \cite{finch2004,naetar2014} and algorithms for spherical and cylindrical measurement setups \cite{haltmeier2005filtered,finch2007inversion,kunyansky2007a} have been developed.
Furthermore, methods that are based on (truncated) series formulas, e.g., eigenfunction expansion, have been developed \cite{agranovsky2007,kunyansky2007b}.
These methods can be implemented with high numerical speed. 
However, they are restricted to certain measurement-surface geometries and typically require additional assumptions, such as constant speed of sound. 
In addition to these methodologies based on analytical inversion formulas, methodologies utilising numerical approximations of the forward model have been developed. These include, for example, time-reversal, regularized least-squares and Bayesian approaches.
Since the wave equation can be stably solved backwards, a physically intuitive approach to image reconstruction is the time-reversal that has been utilised in PAT in various studies, e.g. \cite{xu2004,hristova2008,treeby2010,cox2010b}, including variants based on Neumann series \cite{stefanov2009thermoacoustic,qian2011efficient}.

In limited-view measurement setups, the above described direct reconstruction methods do not perform sufficiently well, and thus variational methods have been utilised. 

In these methods, the reconstructed image is obtained as minimiser of a least squares difference between the measurements and forward model predictions together with a penalising regularisation term \cite{scherzer2009variational}.
In PAT, smoothness supporting  Tikhonov regularisation, the total variation penalty, which supports  piece-wise constant solutions, \cite{arridge2016b,bergounioux2014,wang2012,buehler2011model,javaherian2017} higher order variants \cite{boink2018framework}, and regularisers based on the $\ell_1$-norm promoting sparsity \cite{haltmeier2016compressed} have been utilised. 
Furthermore, the Bayesian approach,  where one aims to solve the posterior distribution of the unknown parameters of primary interest given information from measurements, forward model, and prior model for the unknown parameters, together with their uncertainties \cite{kaipio05}, has been studied, see e.g., \cite{tick2016,tick2018,sahlstrom2020modeling}.
Time-reversal, variational and Bayesian methods are computationally intensive since the wave field in the entire imaging domain needs to be  evaluated frequently.
For more information on approaches for the PAT inverse problem, see e.g. the references above, the reviews \cite{kuchment2011,poudel2019survey,wang2011photoacoustic}, and a review on utilising machine learning methodologies in PAT \cite{hauptmann2020deep}.

The computationally expensive nature of providing PA-images is still one of the bottlenecks of PAT, especially in realistic imaging scenarios that deal with limited sensor geometries and unknown model parameters, as well as in quantitative PAT \cite{cox2012a,tarvainen2024quantitative}. 

In this work we propose an efficient wavelet-based approach for smoothness-promoting reconstructions in PAT. Particular emphasis is placed on exploring the connections between Sobolev penalisation in the $H^s$-scale and smoothness-promoting priors in a Bayesian framework. Beside significantly reducing the memory requirements, the proposed wavelet-based approach significantly reduces the computational workload of the reconstruction process, enabling the computation of large-scale reconstructions.

The remainder of the paper is organized as follows. In Section 2, we formulate a linear mathematical model for the acoustic PAT problem. In Section 3, we give a general overview for solving the inverse problem deterministically, and discuss the relation between specific smoothness promoting regularization methods and the Bayesian framework with Mat{\'e}rn priors. Section 4 discusses  practical implementations for the proposed solution method. In particular, we give a collection of the most relevant definitions and results regarding wavelet transforms, as well as sparse implementation for solving the original inverse problem. In Section 5, the approach is verified numerically.

\section{Photoacoustic tomography}

\subsection{Notations}

Throughout this article, we use the following notations. The space of square integrable functions over a domain $\D \subset \R^d$, satisying $\int_{\D} |u(x)|^2 \d x < +\infty $ is denoted by $L^2(\D)$, with the norm $\|\cdot\|_2$ given as 
$$
\left\| u \right\|_2 = \left( \int_\D |u(x)|^2 \d x\right)^{1/2}.
$$
For all functions $u\in L^2(\R^d)$, the Fourier transform is defined as 
$$
\mathcal{F}u (\xi) = \int_{\R^d} \e^{-\rmi 2\pi \langle \xi,x\rangle} u(x) \d x,
$$
with $\langle \cdot,\cdot \rangle$ being the standard inner-product $\langle x,y\rangle = \sum_{k=1}^d x_k y_k$, in $\R^d$. 
The inverse Fourier transform $\mathcal{F}^{-1}$ maps $u\in L^2(\R^d)$ as
$$
\mathcal{F}^{-1} u(x) = \int_{\R^d} \e^{\rmi 2\pi \langle x,\xi\rangle } u(\xi) \d \xi.
$$
The space of functions satisfying
$$
\mathcal{F}^{-1}\left((1+4\pi^2|\cdot|^2)^{r/2} \mathcal{F}(u)\right)(\cdot) \in L^2(\D),
$$
for all $0\le r\le s$ is called the $L^2$-based Sobolev space with smoothness index $s>0$, and it is denoted by $H^s(\D)$.

\subsection{Forward model}
\begin{figure}
    \centering
    \includegraphics[width=\linewidth]{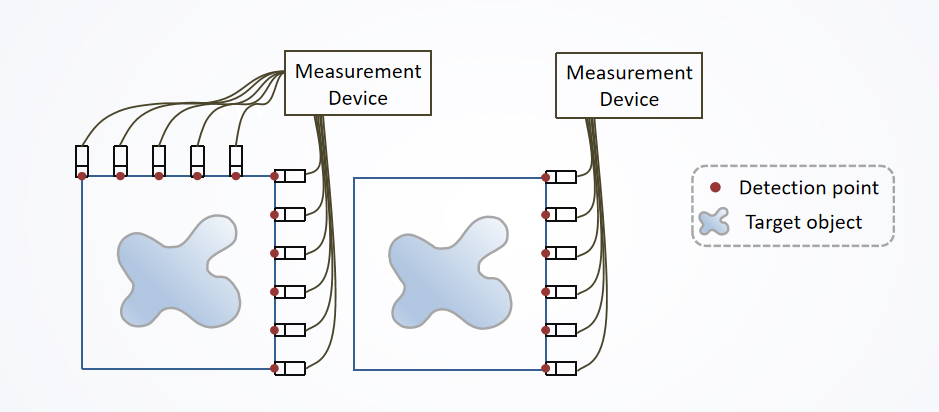}
    \caption{In PAT, the initial pressure of the system relaxes via acoustic wave propagation. The pressure is measured on the boundary of the target over sufficiently long time period, and the task is to reconstruct the initial state of the system (the initial pressure distribution). Due to physical limitations, sensors can not generally cover the entire boundary, leading to PAT problems with limited view. Two- and one-sided geometries for a two dimensional target are depicted here.}
    \label{fig:setups1}
\end{figure}

In PAT experiments, a light pulse is directed into the soft biological tissue under investigation, and the resulting acoustic pulse is measured at the tissue surface, see illustration in Figure \ref{fig:setups1}. 
Assuming an acoustically homogeneous and non-attenuating medium $\D \subset \R^d$, the acoustic wave propagation can be modeled in terms of an initial value problem as
\begin{align}
    \label{eq:wave_equation}
    \begin{cases}
    \nabla^2 p(r,t) - \frac{1}{v^2} \frac{\partial^2 p(r,t)}{\partial t^2} = 0, \\
    p(r, t=0) = p_0(r), \\
    \frac{\partial}{\partial t} p(r, t=0) = 0,
    \end{cases}
\end{align}
where $\nabla^2$ is the Laplace operator, $p(r,t) \in X(\D) \times \mathcal{T}$ is the pressure at point $r\in \D$ and time instance $t\in \mathcal{T}$ ($\mathcal{T} = [0,T]$ for some $T$), $p_0(r)\in X(\D)$ is the initial pressure distribution (IPD), and $v$ is the speed of sound. Here $X(\D)$ denotes a general function space over the domain $\D$, to be specified later.

Let us denote by $\widetilde{A}$ the operator that maps the initial value $p_0$ to the time-dependent solution of the partial differential equation \eqref{eq:wave_equation}. The pressure at the surface of the domain $\partial \D$ is then given in terms of a 'trace' operator 
$$
T: X(\D)\times T \rightarrow X(\partial \D) \times T,
$$
i.e., the restriction of the pressure field to the boundary.

The forward model for full-view acoustic PAT is then obtained as composition
$$
A : = T\circ \widetilde{A}: X(\D) \rightarrow X(\partial \D) \times T, \  p_0 \mapsto T (\widetilde{A} p_0) = p(r_{\partial\D},t).
$$

The focus of this paper is on the limited-view measurement setups. A finite number of point-like ultrasound sensors are placed on some part of the boundary $\gamma \subset \partial \D$. Sensors are located at $r_i \in  \partial \D$, with $i= 1, \ldots , N_s$, and the pressure is measured over finite time period $[0,T]$. 

These sensors record the pressure at finitely many time points, and the results are collected in a matrix form $p_t \in \R^{N_s \times N_t}$, where $N_s$ and $N_t$ denote the number of sensors and the number of time instances, respectively. The physical model is then given in terms of the linear mapping 
$$
K: X(\D) \rightarrow \R^{N_s\times N_t}, p_0 \mapsto [(A p_0)(r_i,t_j)]_{j=1, k=1}^{N_s , N_t},
$$
and the observation are of the form
$$
p_t^\delta = K p_0 + \delta,
$$
where $\delta$ is the measurement noise.

Reconstruction of the IPD $p_0\in L^2(\D)$ can then be obtained as a solution to the linear inverse problem 
$$
K p_0 = p_t^\delta.
$$

To approximate the forward mapping $K$ in a pixel-basis, we use a pseudospectral $k$-space method implemented with the k-Wave MATLAB toolbox, \cite{treeby2010a}.

The full matrix representation of $K \in \R^{N_s N_t \times N}$ can then be compiled by simulating ultrasound waves for each characteristic function $\chi_k$ separately.
While impractical for high-resolution and three-dimensional reconstructions, compiling the full matrix representation of $K$ does allow us to analyse the mapping properties of the forward operator, e.g., for different sensor geometries. To this end, we note that the PAT forward operator is well-conditioned for sufficiently good sensor coverage. Moreover, the inverse problem can be stably solved by some generalized inverse of the forward mapping $K$.

\section{Reconstruction methods}

As mentioned in the previous section, reconstructions for the IPD $p_0$ amounts to the computation of a solution to the linear inverse problem 
\begin{equation}\label{eq:lin_inverse}
    K p_0 = p_t^\delta.
\end{equation}

In the following discussions we introduce the deterministic approach to inverse problems. Therefore we switch to denoting the unknown parameter $p_0 =: x$.

The straightforward approach is to directly solve the normal equation
\begin{equation}
 K^* K x = K^* \pt.
 \label{eq:normal_equation}
 \end{equation}

\begin{figure}
    \centering
    \includegraphics[width=0.6\linewidth]{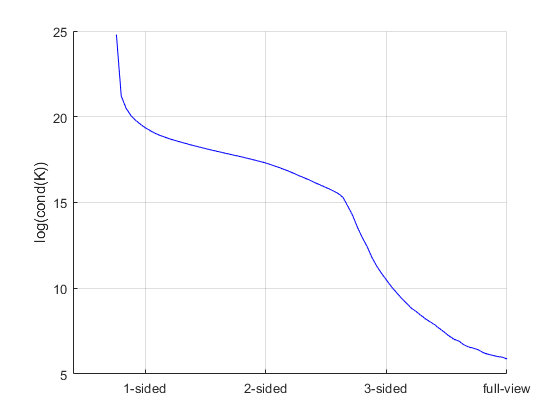}
    \caption{The condition numbers of the acoustic forward operator $K$ vs. the number of sensors at the boundary as described above.}
    \label{fig:condno}
\end{figure}

In general, ill-posed problems can cause the normal equation \eqref{eq:normal_equation} to be unstable to invert, often resulting in issues such as amplification of noise. However, for full-view and near-full-view sensor geometries, the inverse problem \eqref{eq:lin_inverse} is well-posed (see again \cite{xu2004b}), and solving the normal equation still yields reliable reconstructions. The relationship between sensor coverage and the condition number of the forward operator $K$ is illustrated in Figure~\ref{fig:condno}. In a pixel basis the matrix representation of the forward operator $K$ is compiled for several different sensor coverages: Sensors are placed on the boundary of a square domain starting from one sensor placed at a corner and then adding point-like sensors sequentially, until a full-view geometry of 100 equidistantly distributed sensors over the entire boundary is obtained.

We focus on the case of partial measurements. As mentioned, this inverse problem is ill-posed and regularization methods have to be used for a stable approximation of a solution to the inverse problem \eqref{eq:lin_inverse}. This could be either done in a deterministic (see, e.g. \cite{Engl,Ramlau_2015_06,louis1989}) or Bayesian setting \cite{kaipio2006statistical}.

In the deterministic setting, a regularized solution $\xad$ to the generalized solution $\xd$ is computed as
\begin{eqnarray}
    \mathcal{R}_\alpha &:& Y\times \R \to X\\
    \xad & = & \mathcal{R}_\alpha (p^\delta),
\end{eqnarray}
where $\alpha = \alpha (\delta,\pd)$ is the regularization parameter, which has to be chosen such that for $\|\pd -p\|\le \delta$ holds and
\begin{equation}
    \lim_{\delta \to 0}\mathcal{R}_{\alpha(\delta,\pd)} = \xd .
\end{equation}
There is an abundance of regularization methods available, and we refer again to \cite{Engl,Ramlau_2015_06} for an - at least partial - overview.\\
In this paper, we focus on {\it variational} regularization, specifically on Tikhonov regularization, where the generalized solution $\xd$ is approximated by  
\begin{eqnarray}
\xad & = & \arg\min_{x\in X} J_{\alpha,\Omega} (x,\pt)   \\
    J_{\alpha,\Omega} (x) & = &\frac{1}{2} \left\Vert K x - \pt \right\Vert_2^2 + \alpha \Omega(x), 
    \label{eq:reg_functional}
\end{eqnarray}
where $\alpha$ is the regularization parameter, and $\Omega$ is a regularisation functional. The Tikhonov type regularization is considered, as (i) \emph{a priori} information about the object can be easily incorporated in to the reconstruction by a suitable  choice of the regularization functional, (ii) it can also be interpreted in a Bayesian setting and (iii) there are several ways of efficiently computing the minimizers. Standard Tikhonov regularisation refers to the choice of $\Omega = \| \Gamma \cdot \|_2^2$ with some operator $\Gamma$ that can be interpreted as an a prior covariance function in a Bayesian setting.  In our PAT application, we will use the expected smoothness of the solution, measured in Sobolev spaces $H^s$, as prior information, i.e.,

$$\Omega (x) = \|x\|_{H^s}^2.$$

In the following, we will use the Sobolev embedding operators 
\begin{eqnarray}
    E_s&:& H^s (D)\to L_2(D),~D\subset \R^N\\
    x&&\mapsto E_s x = x.
\end{eqnarray}
to characterize and evaluate the $H^s$ norm.
Although formally $E_s$ is the identity (though measured in different metrics), its adjoint $E_s^\ast : L_2(D)\to H^s(D)$ is not. It can be shown that \cite{Ramlau_2004_01,Ramlau_2004_02,HubmerSherinaRamlau}

\begin{eqnarray}
    \|x\|_{H^s} &=& \| \Gamma x \|_{L_2(D)}\\
    \Gamma x &=& \left (E_s^\ast \right )^{-1/2}x
\end{eqnarray}

The choice $\Gamma = (E_s^*)^{-1/2}$ for Tikhonov regularization can be viewed as restricting the domain of the forward operator $K$ to $H^s(\R^{N \times 1}) \subset L^2(\R^{N\times 1})$, which therefore has a regulatization effect, as we are selecting only solutions that belong to $H^s$.

The minimizer $\xad$ of the Tikhonov functional $J_{\|\cdot\|_s^2,\alpha}$ can be computed as solution  to the linear equation
\begin{equation}
(E_s^* K^* K + \alpha I ) x^\dagger = E_s^* K^* \pt ,
\label{eq:linear_minimization}
\end{equation}
where $I$ denotes the identity matrix,  see again \cite{Ramlau_2004_01, HubmerSherinaRamlau}.

An alternative approach toward smoothness promoting reconstructions in  the Bayesian framework has been implemented in \cite{tick2017,tick2018}. Assuming that both the noise and IPD $x$ are mutually independent normally distributed random variables, $x \sim \mathcal{N}(\eta_{x},\Gamma_{x})$ and $\delta \sim \mathcal{N}(\eta_{\delta},\Gamma_{\delta})$,  the reconstructions are obtained as the posterior distribution. 
The posterior is also a Gaussian distribution $x|p_t\sim \mathcal{N}(\eta_{x|p_t},\Gamma_{x|p_t})$ 
and can be written as
\begin{equation}
	\pi(x \vert p_t) \propto \exp \left\{ -\frac{1}{2} \left\Vert L_e(p_t - Kx - \eta_e) \right\Vert^2_2 -\frac{1}{2} \left\Vert L_{x}(x - \eta_{x}) \right\Vert^2_2  \right\}
	\label{eq:Posterior_exp}
\end{equation}
where $L_{e}$ and $L_{x}$ are the Cholesky decompositions of the noise and prior covariance  $\Gamma_{e}^{-1} = L_{e}^{\mathrm{T}}L_{e}$ and  $\Gamma_{x}^{-1} = L_{x}^{\mathrm{T}}L_{x}$.

Reconstructions can now be obtained as the \emph{maximum a posteriori} estimate, i.e., the maximizer $x^\dagger \sim \mathcal{N}(\eta_{x},\Gamma_{x})$ to $\pi(x\vert p_t)$, leading us to the minimization problem

\begin{equation}
x^\dagger = \arg\min_{x\in X} \left(\frac{1}{2} \left\Vert L_e(p_t - Kx - \eta_e) \right\Vert^2_2 +\frac{1}{2} \left\Vert L_{x}(x - \eta_{x}) \right\Vert^2_2  \right).
\label{eq:MAP-estimate}
\end{equation}

For white noise $\delta \sim \mathcal{N}(0,\beta)$, with constant standard deviation 
$\beta$, and prior $x \sim \mathcal{N}(0,\Gamma_{x})$, the minimization in \eqref{eq:MAP-estimate} simplifies to 
\begin{equation} 
x^\dagger = \arg\min_{x\in X} \left(\frac{\beta^2}{2} \left\Vert p_t - Kx \right\Vert^2_2 +\frac{1}{2} \left\Vert L_{x}x \right\Vert^2_2  \right),
\label{eq:MAP-estimate2}
\end{equation}
corresponding to the standard regularization scheme with $\Omega(x) = \|L_x x\|_2^2$ and $\alpha = 1/\beta^{2}$.

A common choice for smoothness promoting prior covariance $\Gamma_{p_0}$ is given in terms of Mat{\'e}rn covariance functions, a class of functions $k_{\nu,\rho}$ defined for positive parameters $\nu$ and $\rho$, and given as
\begin{equation}
k_{\nu,\rho}(\cdot) = \frac{2^{1-\nu}}{\mathbf{G}(\nu)} \left( \frac{\sqrt{2\nu}|\cdot|}{\rho}\right)^\nu K_\nu \left(\frac{\sqrt{2\nu}|\cdot|}{\rho}\right),
\label{eq:matern_covariance}
\end{equation}
where $\mathbf{G}$ denotes the gamma function, and $K_\nu$ is the Macdonald function of order $\nu$.
The corresponding Mat{\'e}rn covariance operator is then given as convolution $\Gamma_{\nu,\rho} : u \mapsto \Gamma_{\nu,\rho} u = k_{\nu,\rho} \ast u $.  
 The relationship between Mat{\'e}rn priors and smoothing on Sobolev scales has been well-established in the literature on Gaussian processes, \cite{rasmussen2006}. 
To further highlight this connection, we give the following proposition:

\begin{proposition}\label{prop}
    For fixed step-size and smoothness parameter $\rho,\nu>0$, the Mat{\'e}rn covariance operator $\Gamma_{\nu,\rho}$ maps $ L^2(\R^d) \rightarrow H^{\nu + \frac{d}{2}}(\R^d)$. 
\end{proposition}

\begin{proof}
    The result follows from the observation that the inverse Fourier transform of the Bessel potential is given as
    \begin{equation}
    \label{eq:bessel_fourier}
    \mathcal{F}^{-1}\left( (1+|\cdot|^2)^{-s/2} \right) (\xi) = \frac{2^{1-\frac{d+s}{2}}}{\pi^{d/2}\mathbf{G}(\frac{s}{2})} K_{\frac{d-s}{2}}(|\xi|) |\xi|^{\frac{s-d}{2}},
    \end{equation}
    see \cite{Aronszajn1961}. It follows that
    \begin{equation}
    k_{\nu,\rho} (x) = \frac{2^d \pi^{d/2}\mathbf{G}(\nu + \frac{d}{2})}{\mathbf{G}(\nu)} \mathcal{F}^{-1}\left( (1+|\cdot|^2)^{-\nu - \frac{d}{2}} \right) \left(\frac{\sqrt{2\nu}}{\rho}x\right),
    \label{eq:Matern_alt}
    \end{equation}
    and therefore, for all functions $u\in L^2(\R^d)$, we have
    $\mathcal{F} \left( k_{\nu,\rho} \ast u\right) = \mathcal{F}k_{\nu,\rho} \cdot \hat{u} \in L_{\nu + d/2}^2(\R^d),$ the weighted Lebesgue space. Thus $k_{\nu,\rho} \ast u \in H^{\nu + d/2}(\R^d)$.
\end{proof}

\begin{proposition}
    With the choice $\rho = \sqrt{2\nu}$, the Mat{\'e}rn covariance operator $\Gamma_{\nu,\sqrt{2\nu}}$ coincides with the adjoint embedding operator $E_{\nu + d/2}^*$ up to a constant multiplier.
\end{proposition}

\begin{proof}
Following \cite{HubmerSherinaRamlau}, for all $u\in L^2(\R^d)$, we have $E_s^*u (x) = \mathcal{F}^{-1}\left((1+|\cdot|^2)^{-s} \hat{u} \right) (x)$ and further, using the formula \eqref{eq:bessel_fourier}, this gives us the characterisation of $E_s^*$ in terms of spatial filter 
\begin{equation*}
    (E_s^* u)(x) = (G_{2s} \ast u)(x),
\end{equation*}
where the filter $G_s$ is given by
\begin{equation}
\label{eq:spatial_filter}
G_s(x) = \frac{2^{1-\frac{d+s}{2}}}{\pi^{d/2}\mathbf{G}(\frac{s}{2})} K_{\frac{d-s}{2}}(|x|) |x|^{\frac{s-d}{2}}.
\end{equation}

In particular, the Ornstein-Uhlenbeck covariance with step-size parameter \mbox{$\rho=1$},
$$
k_{1/2,1}(x) = \frac{\sqrt{2}}{\mathbf{G}(\frac{1}{2})} \left(|x|\right)^{1/2} K_{1/2}\left(|x|\right)(|x|) = \exp{\left(-|x|\right)}
$$
leads to $s=3/2$ smoothing on Sobolev scale, in $\R^2$, cf. 
$$
G_{2\cdot \frac{3}{2}}(x) = \frac{2^{-\frac{3}{2}}}{\pi \mathbf{G}\left(\frac{3}{2}\right)} K_{1/2}(|x|) |x|^{1/2} = \frac{\exp{(-|x|})}{2\pi}.
$$
\end{proof}
\begin{remark}
    In the context of solving the minimisation problems \eqref{eq:linear_minimization} and \eqref{eq:MAP-estimate}, the constant multiplier is absorbed in the regularisation parameter, thus leading to equal outcome.
\end{remark}

\section{Implementation}

In this section, we introduce an efficient and matrix free method for solving the equation \eqref{eq:lin_inverse} in $\R^d$. To estimate the solution, we resort to iterative solution algorithms, e.g., the generalized minimal residual method or the conjugate gradient method. 
As mentioned in the Section 2, the forward operator $K: \R^{N\times 1} \rightarrow \R^{N_s N_t \times 1} $ can be approximated using a $k$-space method, e.g., by using the k-Wave MATLAB toolbox. 
Similarly, the adjoint operator $K^* :\R^{N_s N_t \times 1} \rightarrow \R^{N\times 1} $ can be approximated using the pseudospectral method.
The adjoint embedding $E_s^*$ can be directly implemented by computing the spatial filter $G_s\in \R^{N\times N}$ according to \eqref{eq:spatial_filter}. However, for high reconstruction resolution (large $N\gg 1$), this approach becomes unfeasible due to high memory requirements. 

In the proof of Proposition \ref{prop}, it was noted that the adjoint embedding operator $E_s^*$ can be characterized in the Fourier domain as a multiplication by the potential $(1+4\pi^2|\cdot|^2)^{-s}$. This allows a simple implementation for the discrete problem as
$$
E_s^*X = \mathrm{ifft}\left( B \cdot \mathrm{fft}(X)\right),
$$
where $B\in\R^{n\times n}$ is the pre-computed matrix representation of the potential $(1+4\pi^2|\cdot|^2)^{-s}$, and 
$\mathrm{fft}$ ($\mathrm{ifft}$) denotes the (inverse) fast Fourier transform. 

In the wavelet based approach, the Fourier basis functions $\mathrm{exp}(-\rmi 2 \pi k x)$ are replaced by wavelets. The wavelet basis is constructed by shifted and scaled versions of a function with some desirable properties, such as regularity, oscillation and locality. In particular, the smoothness of a function $f\in L^2(\R^d)$ on a Sobolev scale can be analyzed by the decay rate of the coefficients with respect to scaling factor of the mother wavelet, cf. decay rate of the Fourier coefficients.

In the standard formulation, the function space $L^2(\R^d)$ is decomposed as direct sum 
\begin{equation}
\label{eq:l2_decomp}
L^2(\R^d) = V_0 \bigoplus_{j=1}^\infty W_j, 
\end{equation}
of the low-frequency approximation space $V_0$ and directional high-frequency detail spaces $W_j$. Construction of this decomposition is based on \emph{multiresolution analysis} of the space $L^2(\R^d)$.

To this end, following \cite{Meyer_1993}, we note that for any positive integer $r\ge 1$ there exists compactly supported scaling function $\varphi \in C^r$ and (basic) wavelets $\psi_l \in C^r$, with $l=1,\ldots, 2^d -1$, so that the shifted and scaled versions, 
\begin{align*}
    &\psi_{j,k}^l(x) = 2^{dj/2} \psi_l(2^jx - k),\\
    & \varphi_k(x) = \varphi(x-k),
\end{align*}
for all $k\in \Z^d$ and $j\in \Z$, form an orthonormal basis for the space $L^2(\R^d)$. We use the following indexing for these functions:

Let $\mathcal{M}$ denote the set of vertices of the $d$-dimensional unit cube $[0,1]^d$, excluding the origin.
For each vertex $\varepsilon \in \mathcal{M}$, we identify a unique basic wavelet $\psi_l$, and we denote $\psi^\varepsilon = \psi_l$. Further, let $\Gamma_j$ denote the sequence of lattices  $\Gamma_j = 2^{-j} \field{Z}^d$ for all non-negative integers $j\ge 0$, with $\overline{\cup_{j\in\field{Z}} \Gamma_j} = \R^d$, and let $\Lambda = \left(\cup_{j\in\field{Z}} \Gamma_j \right) \setminus \left\{ 0\right\}$.
The set $\Lambda$ forms indexing of the set of wavelets via the relation $\lambda = 2^{-j} k + 2^{-j-1}\varepsilon \in \Lambda$, i.e.,
\begin{equation}
\label{eq:multid_wavelet}
\psi_\lambda (x) := 2^{dj/2}\psi^\varepsilon (2^j x - k).
\end{equation}
Finally, the following definition gives us the partition of $L^2(\R^d)$ in the form of \eqref{eq:l2_decomp}, and links multiresolution analysis to the wavelet construction from above
\begin{definition}
    Let $\{ V_j \}_{j\in\field{Z}}$ be an r-regular multiresolution analysis of $L^2(\R^d)$ with scaling function $\varphi$, and for all $j\in\field{Z}$, let $W_j$ denote the complement set $W_j = V_{j+1}\setminus V_j$. Moreover, for all $\lambda \in \Gamma_j$ let $\varphi_\lambda (x)= 2^{d j} \varphi(2^j x - k)$. 
    We say that a family of wavelets $\{ \psi_\lambda\}_{\lambda\in\Lambda}$ corresponds to the $r$-regular multiresolution analysis $\{ V_j \}_{j\in\field{Z}}$, if the following conditions hold:
    \begin{enumerate}
        \item the set $\{ \varphi_\lambda\}_{\lambda\in\Gamma_j}$ forms an orthonormal basis of $V_j$, and
        \item the set $\{ \psi_\lambda\}_{\lambda\in\Lambda_j}$ forms an orthonormal basis of $W_j$,
    \end{enumerate}
    where $\Lambda_j = \Gamma_{j+1} \setminus \Gamma_j$.
\end{definition}

It follows that sets $V_0$ and $W_j$ form a partition of $L^2(\R^d)$ of form \eqref{eq:l2_decomp}, i.e., for each $u\in L^2(\R^d)$ there holds 
$$
u  = \sum_{\lambda\in\Gamma_0} \langle u,\varphi_\lambda\rangle_{2} \varphi_\lambda 
+ \sum_{j\ge 0}\sum_{\lambda\in\Lambda_j}   \langle u,\psi_\lambda\rangle_{2} \psi_\lambda.
$$

Now we are in the position to formulate
\begin{proposition} \cite{HubmerSherinaRamlau}
    Let $s\ge 0$, and let the orthogonal wavelet family $\{ \psi_\lambda\}_{\lambda\in\Lambda}$ correspond to an $r$-regular multiresolution analysis of $L^2(\R^d)$ with $s<r$, and with scaling function $\varphi$. Then, for all $u\in L^2(\R^d)$, the following holds:
    \begin{align*}
    E_s^* u &= E_s^* \left(\sum_{\lambda\in\Gamma_0} \langle u,\varphi_\lambda\rangle_{2} \varphi_\lambda 
+ \sum_{j\ge 0}\sum_{\lambda\in\Lambda_j}   \langle u,\psi_\lambda\rangle_{2} \psi_\lambda \right) \\&= \sum_{\lambda\in\Gamma_0} \langle u,\varphi_\lambda\rangle_{2} \varphi_\lambda 
    + \sum_{j\ge 0}\sum_{\lambda\in\Lambda_j} 2^{-2js}  \langle u,\psi_\lambda\rangle_{2} \psi_\lambda
    \end{align*}
    \label{prop:adjoint_wavelet}
\end{proposition}

\begin{corollary}
    Let parameters $\rho, \nu >0 $ be chosen such that $\frac{\sqrt{2\nu}}{\rho} = 2^{-m}$ for some $m\in\field{N} $, and let $u\in L^2(\R^d)$. Then under the assumptions of Proposition \ref{prop:adjoint_wavelet}, the Mat{\'e}rn covariance of the function $u$ agrees to
    \begin{align*}
        \Gamma_{\nu,\rho} u = 2^{d + md}\pi^{d/2} \frac{\mathbf{G}(\nu + \frac{d}{2})}{\mathbf{G}(\nu)} \left( \sum_{\lambda\in\Gamma_{-m}}\langle u,\varphi_\lambda\rangle_{2} \varphi_\lambda +\sum_{j\ge-m} \sum_{\lambda \in \Lambda_j}2^{-2(j+m)(2\nu+d)}  \langle u,\psi_\lambda\rangle_{2} \psi_\lambda\right)
    \end{align*}
\end{corollary}

\begin{proof}
    From equation \eqref{eq:Matern_alt} we have that 
    $$
    k_{\nu,\rho}(\cdot) = 2^d\pi^{d/2} \frac{\mathbf{G}(\nu + \frac{d}{2})}{\mathbf{G}(\nu )} G_{2\nu+\rho}\left(2^{-m} \cdot\right).
    $$
    Now for the convolution $k_{\nu,\rho} \ast u$, we have 
    \begin{align*}
        k_{\nu,\rho} \ast u (x) &= 2^d\pi^{d/2} \frac{\mathbf{G}(\nu + \frac{d}{2})}{\mathbf{G}(\nu )} \left(G_{2\nu+d}\left(2^{-m}\ \cdot\right) \ast u \right) (x) \\
        & = 2^d\pi^{d/2} \frac{\mathbf{G}(\nu + \frac{d}{2})}{\mathbf{G}(\nu )} 2^m \left(G_{2\nu+d} \ast u\left(2^m \cdot\right)\right)\left(2^{-m} x\right) \\
        & = 2^d\pi^{d/2} \frac{\mathbf{G}(\nu + \frac{d}{2})}{\mathbf{G}(\nu )} 2^m \left(E_{\nu+d/2}^* \tilde{u} \right)\left(2^{-m} x\right),
    \end{align*}
    where $\tilde{u}(\cdot) = u\left(2^m \cdot\right)$.

    It follows from Proposition \ref{prop:adjoint_wavelet} that 
    \begin{align*}
    E_{\nu+d/2}^* \tilde{u}  &= \sum_{j<0}\sum_{\lambda\in\Lambda_j} \langle \tilde{u},\psi_\lambda\rangle_{2} \psi_\lambda 
    + \sum_{j\ge 0}\sum_{\lambda\in\Lambda_j} 2^{-j(2\nu + d)}  \langle \tilde{u},\psi_\lambda\rangle_{2} \psi_\lambda,
    \end{align*}
    where the wavelet coefficients are given as
    $$
    \langle \tilde{u},\psi_\lambda\rangle_{2} = 2^{-m}\langle u,\psi_\lambda(2^{-m} \cdot)\rangle_{2} = 2^{-m+dm/2} \langle u,\psi_{\tilde{\lambda}}\rangle_{2}.
    $$
    For a given index $\lambda = 2^{-j}k + 2^{-j - 1}\varepsilon \in \Lambda_j$, the corresponding $\tilde{\lambda} = 2^{-(j-m)}k + 2^{-(j-m) - 1}\varepsilon\in \Lambda_{j-m}$. The claim follows now immediately from the scaling property \eqref{eq:multid_wavelet} of wavelets $\psi_\lambda$.
\end{proof}

In the discrete setting, the wavelet coefficients $\langle u,\varphi_\lambda\rangle_{2}$ and $\langle u,\psi_\lambda\rangle_{2}$ are evaluated by passing the discrete function $X\in L^2(\R^{N})$ repeatedly through low- and high-pass filters. 

Let us consider discretisation of a square $[-a,a]^2\subset \R^2$ into $N = n^2$ pixels. For a fixed wavelet basis, the discrete wavelet transform of a depth $0< m \le \lceil \log_2(n)\rceil$ in this basis is then defined as mapping
$$
\mathcal{W}: L^2(\R^{n\times n}) \rightarrow L^2(\R^{1 \times N_w}), \ X \mapsto \left[ cA_m, cH_m, cV_m, cD_m, cH_{m-1}\ldots, cV_1, cD_1 \right],
$$
where $cA_m$ are the approximation coefficients at depth $m$, and $cH_i, cV_i, cD_i$ are the horizontal, vertical, and diagonal detail coefficients at depth $i=1,\ldots,m$, respectively.

The length of each string $cX_i$ is dependent on the depth $i$. In the MATLAB implementation these lengths are recorded in a book-keeping matrix. We denote the length of string $cX_i$ as $l(i)$.
The discrete wavelet transform depends on the mother wavelet (the wavelet filter). 

For a fixed smoothness parameter $s>0$, and sufficiently smooth mother wavelet ($\psi\in H^r$, where $r>s$), the adjoint embedding operator $\Tilde{E_s^*}$ maps the wavelet coefficients as 
$$
\Tilde{E_s^*} cX_i = 2^{- 2s(m-i)} cX_i.
$$

Therefore, in the wavelet basis, the adjoint embedding operator is a  multiplication with a diagonal matrix with a block structure,
$$
\mathrm{D} = \mathrm{diag}\left(\mathbf{1}_{l(m)},2^{-2s} \mathbf{1}_{3 l(m)}, 2^{-2\cdot 2s}\mathbf{1}_{3 l(m-1)}, 2^{-3\cdot 2s}\mathbf{1}_{3l(m-2)},\ldots, 2^{-(m-1)\cdot 2s}\mathbf{1}_{3l(1)}\right),
$$
where $\mathbf{1}_n = (1,1,\ldots, 1)$ of length $n$.
Therefore the original adjoint embedding operator $E_s^* : l^2(\R^{n\times n}) \rightarrow H^s(\R^{n\times n})$ can be efficiently implemented as 
$$
E_s^* X = \mathcal{W}^{-1} \left( \mathrm{D} \cdot (\mathcal{W} X) \right),
$$
where both $\mathcal{W}$ and $\mathcal{W}^{-1}$ are implemented via fast wavelet transform.

\section{Simulations}

In order to simplify the interpretation of the sensor coverage, a  $50 \, \mathrm{mm} \times 50 \, \mathrm{mm}$ square domain was considered with two sensor geometries where sensors were located on one side and on two adjacent sides of the object. The speed of sound was set to $v = 1500 \, \mathrm{m/s}$. These setups are illustrated in Figure \ref{fig:setups1}, with total of 80 equidistantly located sensors in both examples.

\begin{figure}[tbp!] 
\centering 
\includegraphics[width=.5\linewidth]{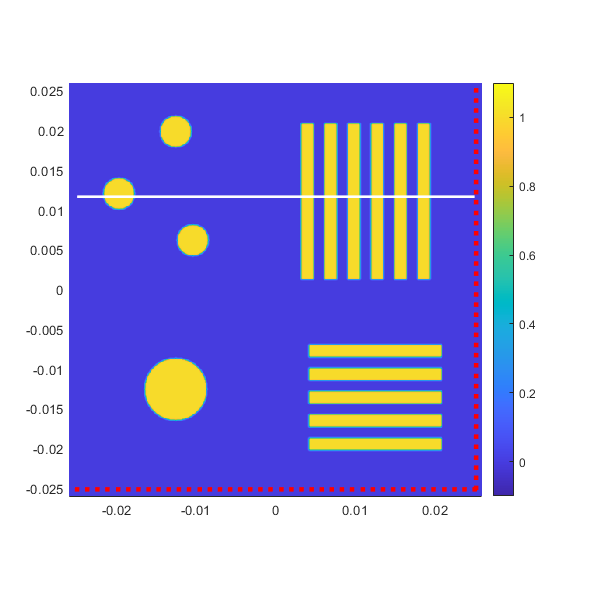} 
\caption{Simulated initial pressure distribution 
$p_0$, consisting of several rectangular and circular inclusions within a square-shaped domain. Sensors are evenly distributed along one and two sides of the target. The white line indicates the cross-section for reconstructions shown in Figure \ref{fig:slice}.} 
\label{fig:setup} 
\end{figure}

The simulated initial pressure distribution $p_0$ consisted of multiple circular and rectangular inclusions within the domain. The initial pressure 
$p_0$ and sensor geometries are visualized in Figure \ref{fig:setup}.

For simulating the data, the spatial domain was discretized into $1000 \times 1000$ pixels with a pixel size of $\Delta x = 50.0 \, \mu\mathrm{m}$. A perfectly matched layer (PML) of 10 pixels was added to absorb boundary reflections. The number of time points was set to 
$N_t = 4809$ with a time step of $10.0 \, \mathrm{ns}$, leading to a total of 384720 data points. Data were simulated using the pseudospectral method implemented in the k-Wave toolbox \cite{treeby2010a}. To model noise, Gaussian white noise with zero mean and a standard deviation of 5\% of the maximum signal amplitude was added.

For the reconstructions, the domain was discretized into $512 \times 512$ pixels with a pixel size of $\Delta x = 97.7 \, \mu\mathrm{m}$. The PML thickness was 10 pixels. The temporal grid was defined with $N_t = 2508$ time points and a time step of $19.5 \, \mathrm{ns}$. The simulated data were linearly interpolated to match the temporal discretization used for reconstructions.

\begin{figure}
    \centering
    \includegraphics[width=1\linewidth]{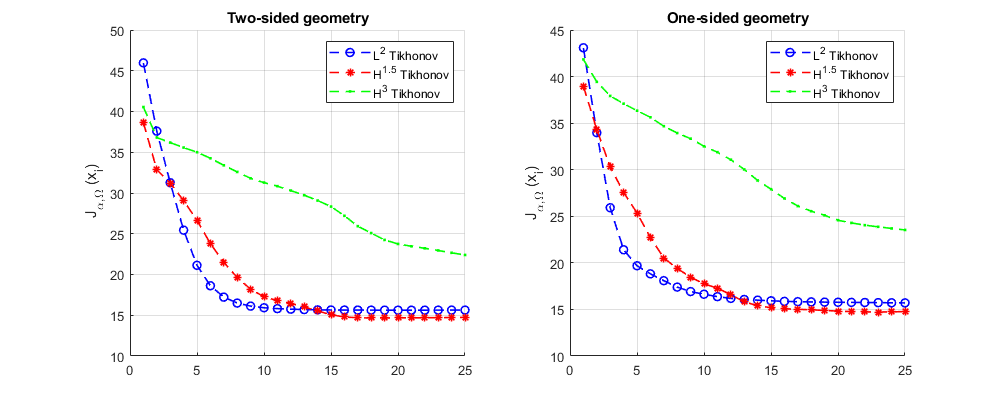}
    \caption{The functional $J_{\alpha,\Omega}$ is minimized using the generalized minimal residual method with smoothness parameters $s=0,1.5,3$. The values of $J_{\alpha,\Omega}$ plotted for each iterate corresponding to minimal residual solution to equation \eqref{eq:linear_minimization}, approximated in Krylov space $\mathcal{K}_i(E_s^* K^* K + \alpha I,  E_s^* K^* \pt )$. The regularization parameter was set at $\alpha=10^{-5} $. }
    \label{fig:func_convergence}
\end{figure}

The solution to equation \eqref{eq:linear_minimization} was estimated using the generalized minimal residual (GMRES) method with 15 iterations, cf. Figure \eqref{fig:func_convergence}, resulting in a total of 31 evaluations of the forward and adjoint wave propagation problems. Reconstructions were computed using smoothness parameters $s=0$ (corresponding to standard $L^2$ Tikhonov regularization), 
$s=3/2$ (corresponding to an Ornstein–Uhlenbeck (O-U) covariance prior), and $s=3$. For reference, total variation (TV) reconstructions were computed using 50 Barzilai–Borwein iterations, requiring 101 evaluations of the wave propagation.

\begin{figure} 
\centering 
\includegraphics[width=.85\linewidth]{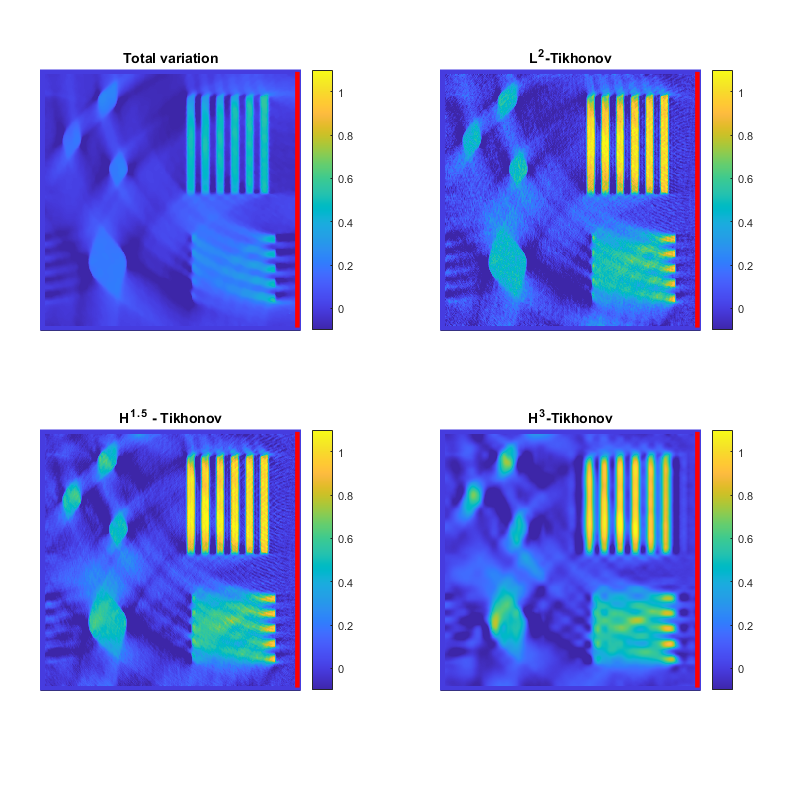} \caption{Estimated $p_0$  in the one-sided sensor geometry. Top row, left to right: TV, $s = 0$, $s = 3/2$ (corresponding to Mat{\'e}rn covariance $\Gamma_{1/2,1}$ in $\R^2$), and $s = 3$ (Mat{\'e}rn covariance $\Gamma_{2,1}$ in $\R^2$).} \label{fig:one-sided}
\end{figure}

\begin{figure}
\centering 
\includegraphics[width=.85\linewidth]{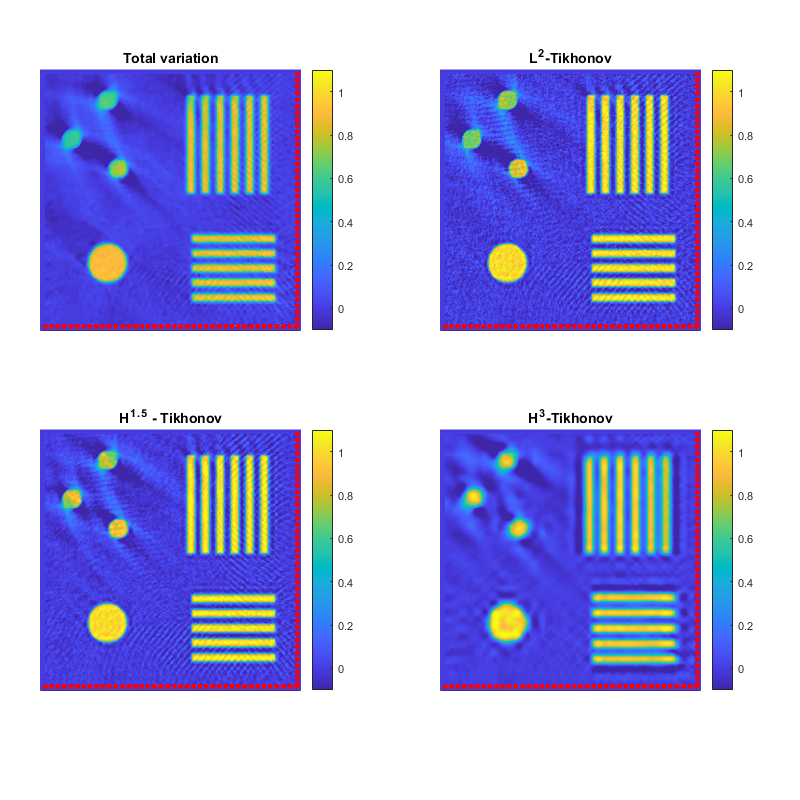} \caption{Estimated $p_0$ in the two-sided sensor geometry. Top row, left to right: TV, $s = 0$, $s = 3/2$ (corresponding to Mat{\'e}rn covariance $\Gamma_{1/2,1}$ in $\R^2$), and $s = 3$ (Mat{\'e}rn covariance $\Gamma_{2,1}$ in $\R^2$).} \label{fig:two-sided} 
\end{figure}

Figures \ref{fig:one-sided} and \ref{fig:two-sided} show the estimated $p_0$ for the one-sided and two-sided sensor geometries, respectively, for TV and the three smoothness parameters. The corresponding relative errors (RE) are listed in Table \ref{tab:RE}, computed after interpolating the estimated parameters to the simulation discretization.

\begin{table}[tbp!] 
\centering 
\caption{Relative error (RE \%) for the estimated $p_0$ using TV and smoothness parameters $s = 0$, $s = 3/2$ and $s = 3$ in one- and two-sided sensor geometries. 
} 
\vspace{.2cm} 
\begin{tabular}{c|c|c|c|c} \hline  
& TV & $s = 0$ & $s = 3/2 $ & $s = 3$  \\ \hline
One side & 0.99 & 0.90 & 0.80 & 0.93 \\
Two side & 0.45 & 0.42 & 0.37 & 0.54 \\ \hline 
\end{tabular} \label{tab:RE} 
\end{table}

The results show that with TV regularization, significant limited view artifacts are present, particularly for the one-sided sensor geometry. In contrast, the estimates with 
$s=0$ and $s=3/2$ display considerably reduced artifacts, both in the background and the sharpness of the inclusions. Notably, reconstructions with $s=3$ yield background values closer to the ground truth compared to TV, but at the cost of a slight smoothing of the inclusions, which is more pronounced for thin, high-contrast structures such as the rectangular inclusions.

\begin{figure} 
\centering 
\includegraphics[width=1\linewidth]{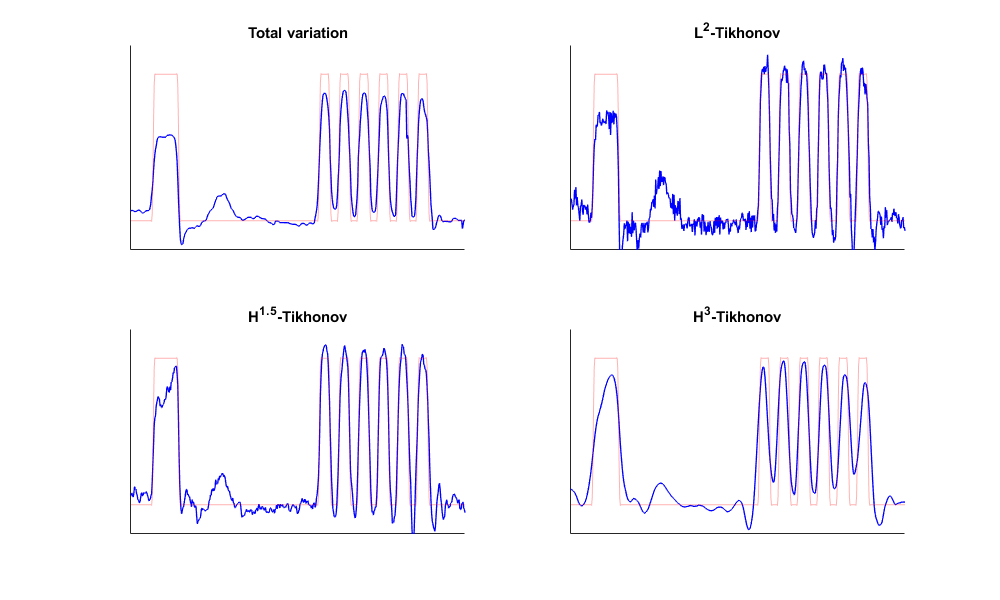} \caption{Cross-section along the white line in Figure \ref{fig:setup}. The red line shows the ground truth and the blue lines show the corresponding reconstructions for TV and $s = 0$, $s = 3/2$ and $s = 3$.}
\label{fig:slice} 
\end{figure}

The cross-sections of the reconstructions for the two-sided sensor case, presented in Figure \ref{fig:slice}, illustrate these effects clearly. Reconstructions using $s=3$ display visibly smoothed peaks compared to the ground truth, particularly for the thinner inclusions.

From the relative errors, it is evident that both standard 
$L^2$- Tikhonov and $H^{3/2}$ smoothing priors provide comparable reconstruction quality, while $H^3$ prior leads to a deterioration due to over smoothing in relation to the known ground truth. This effect becomes more pronounced with increasing distance from the sensors.

It is also worth noting that TV regularization tends to compress contrast, leading to underestimated peak values of the inclusions. Moreover, the computational efficiency of the proposed method is demonstrated by the number of wave equation evaluations: only 31 evaluations were needed for our approach compared to 101 for TV reconstructions, reflecting a significant reduction in computational cost. 

\section{Discussion}

In our paper the inverse problem of PAT was considered. A special emphasis was given to linking Bayesian approaches to well-established results stemming from regularization based methods. In particular, the connection between the \emph{maximum a posteriori} estimate and Tikhonov regularized solutions has been discussed. In both occasions an optimization problem has to be solved, and its solution is approximated with iterative algorithms. Additionally, an efficient wavelet based implementation of the smoothing Mat{\'e}rn covariance operators was established, which is based on adjoint Sobolev embedding operators expressed in terms of weighted wavelet coefficients.

In the proposed method the implementation of the Sobolev / Mat{\'e}rn prior is done at a computational cost of $O(N_x)$, where $N_x$ is the number of pixels in the reconstruction. 
The linear dependence on the number of pixels allowed us to upscale the problem size 16-fold, from an image size of $128 \times 128$ pixels to $512 \times 512$ pixels, when computations were performed on regular laptop. In contrast this was impossible due to memory limitations as the pre-computed Ornstein-Uhlenbeck covariance in the pixel basis scales quadratically $O(N_x^2)$ with the number of pixels. The limit in RAM was reached with discretisation of $170\times 170$, with covariance requiring 3.34GB when compiled as a single-precision matrix.

Numerical examples included reconstructions of the initial pressure distribution for several different smoothness indices for the solution and two different sensor geometries. The reconstructions were obtained by iteratively estimating the solution to the equation \eqref{eq:linear_minimization}. Both the forward and adjoint operators were implemented in a matrix-free manner, allowing for higher discretisations to be considered. The used simulation environment for the acoustic wave propagation required a pixel base to be used, and thus the covariance had to be implemented using fast wavelet transform at each iteration. In the context of general inverse problems, additional improvements in both memory- and computational efficiency could be reached by solving the problem directly in a wavelet basis, e.g., by using a wavelet Galerkin method. However, further numerical validation at larger scales would be needed confirm this.

\section{Acknowledgments}
Authors JK and RR were funded in part by the Austrian Science Fund (FWF) SFB \url{https://doi.org/10.55776/F68} ‘Tomography Across the Scales’, project F6805-N36 (Tomography in Astronomy). For open access purposes, the authors have applied a CC BY public copyright license to any author-accepted manuscript version arising from this submission.

Authors TS and TT have also received funding from the European Research Council (ERC) under the European Union’s Horizon 2020 research and innovation programme (grant agreement No 101001417 - QUANTOM) and the Research Council of Finland (Centre of Excellence in Inverse Modelling and Imaging grant 353086, Flagship of Advanced Mathematics for Sensing Imaging and Modelling grant 358944, and Flagship Program Photonics Research and Innovation grant 320166).

\bibliographystyle{plain}
\bibliography{pat}

\end{document}